\documentclass[12pt]{hha}
\usepackage{diagrams,amssymb,emlines,sphat}
\usepackage{amsfonts,amsthm,amsmath}
\begin{document}{}
\hfuzz=30pt
\title{Group Extensions And Automorphism Group Rings\\}


\date{\today }

\author{John Martino}
\email{martino@wmich.edu}       
\address{Department of Mathematics\\ 
         Western Michigan University\\
         Kalamazoo, MI 49008\\
         U.S.A.}

\author{Stewart Priddy}            
\email{priddy@math.northwestern.edu}       
\address{Department of Matheamtics\\        
         Northwestern University\\
         Evanston, IL 60208\\
         U.S.A.}



\begin{abstract}
We use extensions to study the semi-simple quotient of the group ring $\mathbf{F}_pAut(P)$
of a finite $p$-group $P$.
For an extension $E: N \to P \to Q$,
our results involve relations between $Aut(N)$, $Aut(P)$, $Aut(Q)$ and the extension class
$[E]\in H^2(Q, ZN)$. One novel feature is the use of the {\it intersection orbit group}
$\Omega([E])$, defined as the intersection of the orbits $Aut(N)\cdot[E]$ and $Aut(Q)\cdot [E]$
in $H^2(Q,ZN)$. This group is useful in computing $|Aut(P)|$.
In case $N$, $Q$ are elementary Abelian $2$-groups our results involve the theory of quadratic forms and
the Arf invariant.
\end{abstract}

\maketitle

\newtheorem{Thm}{Theorem}[section]
\newtheorem{Lem}[Thm]{Lemma}
\newtheorem{Cor}[Thm]{Corollary}
\newtheorem{Prop}[Thm]{Proposition}
\newtheorem{Ex}[Thm]{Example}
\newtheorem{Def}[Thm]{Definition}
\theoremstyle{definition}
\newtheorem{Rem}[Thm]{Remark}

\thanks{}



\section{Introduction}

Since the simple modules of a ring and its semi-simple quotient are the same, for many
purposes it suffices to consider the latter ring. In this note we study 
the problem of calculating the semi-simple quotient of the group ring
$\mathbf F_p Aut(P)$ for the automorphism group of a finite $p$-group. The usual method is to 
consider the maximal elementary Abelian quotient
$P \to P/\Phi(P)$, where $\Phi(P)$ is the Frattini subgroup.
The induced map $Aut(P) \to Aut(P/\Phi(P))$ has a $p$-group as kernel by the Hall Basis
Theorem. Hence the map of algebras
$$       
         \mathbf{F}_pAut(P) \to \mathbf{F}_pAut(P/\Phi(P))  
$$
has a nilpotent kernel and thus suffices to compute the  semi-simple quotient. However this map
is not necessarily onto and one is left with the still considerable problem of determining
the image. For $\Phi(P) = \mathbf{Z}/p$, when
$p$ is an odd prime,
this has been done by Dietz \cite{D} giving a complete determination of $Aut(P)$.

We adopt an inductive approach via extensions; that is, we assume $P$ is given as an extension
$$
     E: N \to P \to Q
$$
with $Aut(N)$, $Aut(Q)$ under control. Then there is an exact sequence relating the 
automorphism groups of $N$ and $Q$ with that of $P$, depending on the cohomology class of the 
extension  $[E] \in H^2(Q,Z(N))$ where $Z(N)$ denotes the center of $N$.

Our motivation comes from stable homotopy theory.
Let $G$ be a finite group and $p$ be a prime.
The classifying space of $G$ completed at $p$, $BG^{\wedge}_p$, decomposes stably into a wedge of 
indecomposable summands
$$ BG^{\wedge}_p \simeq X_1 \vee X_2 \vee \cdots \vee X_n.$$
Each summand $X_i$ is the mapping telescope of a primitive idempotent $e$ in the ring of stable self-maps 
of $BG^{\wedge}_p$, $e\in\{ BG^{\wedge}_p,BG^{\wedge}_p\}$.  Thus there is a one-to-one correspondence 
between the indecomposable summands and the
simple modules of the ring of stable self-maps.
This correspondence is explored in both \cite{BF} and \cite{MP} (see also \cite{MP1}). It turns out that
modular representation theory plays a crucial role:
if $P$ is a Sylow $p$-subgroup of $G$ 
then each indecomposable summand of $BG^{\wedge}_p$ originates in 
$BQ$ for some subgroup $Q \leq P$ and corresponds to a simple $\mathbf{F}_pAut(Q)$ module.

Of course, the automorphism group of a group is of intrinsic interest in its own right, and our 
methods shed some light on its structure.

An outline of the paper follows: Section 2 covers the preliminaries on extensions,
$E: N \to G \to Q$,
including the fundamental exact sequence, Theorem \ref{ThmWB}, relating
$Aut(N)$, $Aut(G)$, $Aut(Q)$, and the extension class $[E]$. 
In Section 3 we define and identify a 
group structure, $\Omega([E])$, on the
intersection of the two orbits $Aut(N)\cdot [E]$ and $Aut(Q)\cdot [E]$ where 
$Aut(N)$ and $Aut(Q)$ act on  $H^2(Q,Z(N))$ in the usual way. This group, which 
we call the intersection orbit group, is useful in computing $|Aut(G)|$.
The case of trivial action (or twisting) of $Q$ on $Z(N)$ is considered in Section 4. 
Extensions with $N$, $Q$, elementary Abelian $p$-groups are studied 
in Section 5. In case $p=2$ this involves the theory of quadratic forms over $\mathbf{F}_2$
and the Arf invariant. We recall Browder's classification theorem  \cite{Br} and give
several results describing the order of a quotient of $Aut(G)$ by a normal $p$-subgroup in Section 6.
For more complicated $p$-groups $G$ we describe an inductive procedure for extending these results
using the mod-$p$ lower central series. Section 7 is devoted to several applications of the theory.

In what follows all groups are assumed finite, except as noted in Section 3.

\medskip

\section{Preliminaries}

We begin by recalling the results of C. Wells \cite{W} as extended by J. Buckley \cite{Bk}.
Because their notation is now non-standard, e.g., functions written on the right,
we re-couch these results in more standard notation.
Let 
$$
         E: N \overset{i}\rightarrow G \overset{\pi}\rightarrow Q
$$
be an extension of the group $N$ by the group $Q$ and let $Aut_N(G)$ be the group of 
automorphisms of $G$ mapping $N$ to itself. The obvious homomorphism
$ \rho = (\rho_Q, \rho_N): Aut_N(G) \to Aut(Q) \times Aut(N)$ provides a means of studying 
$Aut_N(G)$.

As usual two extensions $E_1$, $E_2$
are {\it{equivalent}} $E_1 \sim E_2$ if there is an isomorphism $\alpha :G_1 \to G_2$ 
restricting to the identity on $N$ and inducing the identity on $Q$. The set of such equivalent 
extensions is denoted $\mathcal{E}(Q,N)$.
The {\it{twisting}} (or {\it coupling}) $\chi :Q \to Out(N)$
of $E$ is the homomorphism defined as usual by $\chi(q)(n) = i^{-1}(g^{-1}i(n)g)$ where $g \in \pi^{-1}(q)$, $n\in N$.
Equivalent extensions have the same twisting. 

The center $ZN$ of $N$ has the structure of a $Q$ module via a homomorphism
$\overline{\chi} : Q \to Aut(ZN)$ defined by the composite
$$
  \overline{\chi}:Q \overset{\chi}\rightarrow Out(N) \overset{res}\rightarrow Aut(ZN)
$$
where $res  : Out(N) \to Aut(ZN)$ is induced by $Aut(N) \to Aut(ZN)$.
It is well-known that we may identify
$$
\mathcal{E}(Q,N) =  \coprod _{\chi} H_{\overline{\chi}}^2(Q, ZN)
$$
where $\chi$ ranges over all  twistings $\{ Q \to Out(N)\}$.

Now consider $(\sigma, \tau) \in Aut(Q) \times Aut(N)$ and form the extension
$$
    \sigma {E} {\tau}^{-1}: N \overset{i{\tau}^{-1}}\rightarrow  G \overset{\sigma \pi}\rightarrow Q
$$
Then $ Aut(Q) \times Aut(N)$ acts on $\mathcal{E}(Q,N)$ from the left
by
\begin{equation}
   (\sigma, \tau)[E] = [\sigma E \tau^{-1}].
\label{eq:a}
\end{equation}

One checks $(\sigma, \tau)(\sigma', \tau') [E] = (\sigma\sigma',\tau\tau')[E]$ and
$(1,1)[E]= [E]$. 
The twisting of $(\sigma, \tau)E$ is given by $\gamma_{\overline{\tau}}\chi\sigma^{-1}$ where 
$\gamma_{\overline{\tau}}$ denotes conjugation by $\overline{\tau}$, the image of 
$\tau$ in $Out(N)$. For a given
$\chi$ define the subgroup $C_{\chi} \subset Aut(Q) \times Aut(N)$ by
$$
C_{\chi} = \{ (\sigma, \tau)\in Aut(Q) \times Aut(N)\ |\ \gamma_{\overline{\tau}}\chi\sigma^{-1} 
= \chi \}
$$
that is, the following diagram commutes
\label{eq: b}
\[                  
\begin{diagram} \label{eqno: 2}
Q &\rTo{\chi} & Out(N)  \\                          
\dTo{\sigma} & & \dTo{\gamma_{\overline{\tau}}} \\
Q & \rTo{\chi} & Out(N)
\end{diagram}
\]
The subgroup $C_{\chi}$ consists of all ordered pairs $(\sigma,\tau) \in Aut(Q) \times Aut(N)$ that preserve the twisting.

If $\chi$ is trivial then clearly $C_{\chi} = Aut(Q) \times Aut(N)$. We note that
$ker(\chi)$ plays no role in the commutativity of the diagram and $\sigma|_{ker(\chi)} 
:ker(\chi) \to ker(\chi)$. Thus if the sequence
$$
      ker(\chi) \to Q \to im(\chi)
$$ 
splits, e.g., if $Q$ is elementary Abelian, then $\sigma|_{ker(\chi)}$ can be
an arbitrary linear isomorphism.

Then $C_{\chi}$ acts on $\{[E]\ |\ E \ has \ twisting \ \chi\}$. 
It is trivial to check that $Im(\rho) \subset C_{\chi}$ so we consider $\rho$ as a 
homomorphism $\rho : Aut_N(G) \to C_{\chi}$.

Let $Z_{\overline{\chi}}^1(Q, ZN)$ denote the group of derivations, i.e., functions
$f: Q \to ZN$ satisfying $f(qq') = 
f(q) + qf(q')$ for $q,q' \in Q$.
Then there is a homomorphism $\mu: Z_{\overline{\chi}}^1(Q, ZN) \to Aut_N(G)$ defined by
$\mu(f)(g) = f(\pi(g))\cdot g$. 

Finally we define a function 
$\epsilon: C_{\chi} \to H_{\overline{\chi}}^2(Q, ZN)$ by restricting the action of 
$Aut(Q) \times Aut(N)$ 
on the extension class $[E]$, that is, $(\sigma, \tau) \mapsto (\sigma, \tau)[E]$.
In general $\epsilon$ is not a homomorphism.

The following is the principal result of \cite{W} as extended by \cite{Bk}:

\medskip
\begin{Thm}
\label{ThmWB} For a given extension $E: N \to G \to Q$ with twisting \linebreak $\chi: Q \to Out(N)$
there is an exact sequence
$$
   1 \to Z_{\overline{\chi}}^1(Q, ZN)\overset{\mu}\rightarrow Aut_N(G) \overset{\rho}\rightarrow
C_{\chi} \overset{\epsilon}\rightarrow H_{\overline{\chi}}^2(Q, ZN)
$$
with $Im(\rho) = (C_{\chi})_{[E]}$,
the isotropy subgroup of $C_{\chi}$ fixing $[E]$. The map $\epsilon$ is not onto and is only a 
set map.

\end{Thm}
\smallskip

An alternate exact sequence results by replacing $Z_{\overline{\chi}}^1(Q, ZN)$ with the cohomology 
group $H_{\overline{\chi}}^1(Q, ZN)$ and $Aut_N(G)$ by $Aut_N(G)/Inn_{ZN}(G)$
where $Inn_{ZN}(G)$ is group of inner automorphisms of $G$ induced by elements of $ZN$, (see \cite{Mac},
Proposition IV 2.1).
We shall be interested in $Im(\rho)$ so we will not need this refinement.

\medskip

In case $ZN$ is a $p$-group Theorem \ref{ThmWB} shows $Z_{\overline{\chi}}^1(Q,ZN)$ is a normal
$p$-subgroup of $Aut_N(G)$ which in turn implies $ \mathbf{F}_p((C_{\chi})_{[E]})$ suffices
to compute simple modules and idempotents:

\begin{Cor}
\label{CorB} If  $E: N \to  G \to Q$ is an extension with $ZN$ a p-group, then
$\mathbf{F}_pAut_N(G)  \overset{\mathbf{F}_p(\rho)}\longrightarrow
\mathbf{F}_p((C_{\chi})_{[E]})$ is surjective with nilpotent kernel.
\end{Cor}
{\it Proof:} Since $ker(\rho) = Z_{\overline{\chi}}^1(Q, ZN)$ is a $p$-group it is well known
that $ker(\mathbf{F}_p(\rho))$ is nilpotent. \qed

\medskip

We shall be interested in determining when $Aut_N(G)$ is a $p$-group; clearly we have

\begin{Cor}
\label{CorBAA} If $Aut(Q)$
and $Aut(N)$ are $p$-groups then $Aut_N(G)$ is a $p$-group.
\end{Cor}

\medskip

\begin{Cor}
\label{CorBA}
Suppose $G$ is a $p$-group with $p$ odd, and $N$ is generated by the elements of $G$ of order $p$.
Then $Aut(G)$ is a $p$-group if $Aut(N)$ is a $p$-group.
\end{Cor}
{\it Proof:} By hypothesis $N=\Omega_1(G)$ is a characteristic subgroup, hence $Aut_N(G) = Aut(G)$.
If $\alpha \in Aut(G)$ has $p'$ order then $\rho_N(\alpha) = id_N$.
Since 
$p$ is odd, the automorphisms of $G$ which have $p'$ order and fix 
$N$ are trivial by Theorem 5.3.10, \cite{Go}. Thus $Aut(G)$ is a $p$-group. \qed

\bigskip

\section{The Intersection Orbit Group}
\smallskip

Let $X$ be a left $A\times B$ set where $A$, $B$ are groups, not necessarily finite. 
Equivalently $A$ acts
on the left, $B$ acts on the right such that $a(xb)=(ax)b$, and $(a,b)x = axb^{-1}$,
where $a\in A$, $x \in X$, $b\in B$. 
As usual let $Ax$
denote the orbit of $x$ under the action of $A$ (respectively $xB$ denote the orbit of
$x$ under the action of $B$). We define the {\it intersection orbit group at $x$}
$$
     \Omega(x) := (Ax) \cap (xB)
$$

If $ax = xb$, $a'x = xb'$ are elements of $\Omega(x)$, their product is defined by
$$
          (ax)(a'x) := (aa')x = x(bb')
$$
It is straightforward to check that this pairing is well-defined giving $\Omega(x)$
the structure of a group. Although left and right actions commute, $\Omega(x)$ is not 
necessarily Abelian; however, if either
$A$ or $B$ is a $p$-group, then $\Omega(x)$ is a $p$-group.
Also if $A$ or $B$ is trivial then obviously $\Omega(x) = \{x\}$,
the trivial group. 

Again, as usual, let  $A_{x}$ and
$B_{x}$ denote the respective isotropy subgroups.
\medskip

\begin{Prop}
\label{PropAZ} There is an isomorphism of groups
$$
\phi :(A\times B)_{x}/(A_{x}
\times B_{x}) \overset{\cong}\longrightarrow \Omega(x)
$$
given by $\phi(a,b) = ax = x{b}^{-1}$.
\end{Prop}
{\it{ Proof:}}
$(a,b)x = x$ if and only if $ax = x{b}^{-1}$.
Similarly $(a,b) \in A_{x}\times B_{x}$ if and only if
$ax = x{b}^{-1} = x$. Thus $\phi$ is well-defined
and bijective. It is a homomorphism by the definition of
the product in $\Omega(x)$. \qed
\medskip

\begin{Cor}
\label{CorAB} If $A$ and $B$ are finite groups, then
\smallskip


$|\Omega(x)|$ divides
$$
|N_{A}(A_{x})/ A_{x}|,
|N_{B}(B_{x})/B_{x}|,
$$
and
$$
 gcd( |Ax|,
|xB|).
$$
\end{Cor}

{\it{Proof:}} Note
that $N_{A}(A_{x})/ A_{x}$ is the largest subset of $Ax \cong A/A_{x}$
which is a group. Thus $\Omega(x) \leq N_{A}(A_{x})/ A_{x}$ is a subgroup
and so $|\Omega(x)|$ divides $|N_{A}(A_{x})/ A_{x}|$ and hence $A/A_{x}$.
Similarly for $B$. Thus
$|\Omega(x)|$ divides $gcd(|Ax|,|xB|)$. \qed
\medskip

Returning to automorphism groups, the following result often simplifies
computing $|Im(\rho)|$.

\begin{Cor}
\label{CorABB} 
For a given extension of finite groups $E: N\to G \to Q$, 

i)  $|Im(\rho)| = |{Aut(Q)}_{[E]}|\cdot |{Aut(N)}_{[E]}| \cdot |\Omega([E])|$

ii) $Im(\rho)$ is a $p$-group
if and only if ${Aut(Q)}_{[E]}$, ${Aut(N)}_{[E]}$, and $|\Omega([E])|$ are $p$-groups.
\end{Cor}

{\it{Proof:}} ii) follows from i) which follows immediately from Proposition \ref{PropAZ}
with $A = Aut(Q)$, $B = Aut(N)$. \qed
\bigskip

There exists extensions that do not depend on the intersection orbit group.  
\begin{Prop}
\label{PropO}
Suppose $N$ is an Abelian $p$-group, $Aut(Q)$ is a $p$-group, and 
$H^2(Q,N) \neq 0$ with a trivial twisting.  Then there exists
a non-split extension $[E]$ such that
$$
|Aut_N(G)|_p = |Hom(Q,N)| |Aut(Q)| |Aut(N)|_p
$$
\end{Prop}

{\it{Proof:}} 
By Corollary \ref{CorAB}, it suffices find $[E]$ with  $|\Omega([E])| =1$.
Let $Aut(Q) \times S$ be a Sylow $p$-subgroup of $Aut(Q) \times Aut(N)$. 
The $Aut(Q) \times Aut(N)$-module $H^2(Q,N)$ has cardinality a power of $p$. Thus there is 
an $Aut(Q) \times S$ fixed point, $[E] \neq 0$ by the fixed point result \cite{Se}, p.\ 64, which generalizes to this case). Thus $Aut(Q)_{[E]} = \{[E]\}$ and $|S| = |Aut(N)|_p$.  The kernel of the exact sequence in Theorem \ref{ThmWB} is $Z^1(Q,N)=Hom(Q,N)$, which explains the presence of that term in the proposition.  \qed

\bigskip

\section{Trivial twisting}

\medskip

\begin{Prop}
\label{PropA}
If $E: N \to G \to Q$
has trivial twisting the exact sequence of Theorem \ref{ThmWB} reduces to
$$
   0\to  Hom(Q, ZN) \overset{\mu}\rightarrow Aut_N(G) \overset{\rho}\rightarrow Aut(Q) 
\times Aut(N) \to H^2(Q, ZN)
$$
\end{Prop} 
{\it Proof:} \cite{Bk}, Th. 3.1. Clearly $Z_{\overline{\chi}}^1(Q, ZN) 
= Hom(Q, ZN)$ and $C_{\chi} = Aut(Q) \times Aut(N)$. \qed
\medskip

Before proceeding, it is instructive to consider the simplest case $G = N \times Q$.
The extension class is trivial thus $\rho: Aut_N(G) \to Aut(Q) \times Aut(N)$
is surjective. In fact
$\rho$ is split by the usual inclusion $Aut(Q) \times Aut(N) \to Aut_N(G)$.
For a more complete discussion of $Aut(Q \times N)$ see \cite{MPD}.
\medskip

In what follows we are interested in $G = P$, a $p$-group. As  an immediate corollary
of  \ref{PropA} we have
\medskip



\begin{Cor}
\label{CorCB} Let $P$ be a $p$-group and define $P_1 = P$, $P_{i+1} = P_i/Z(P_i)$. 
Suppose $Aut(Z(P_i))$ is a $p$-group for each $i$. Then $Aut(P)$ is a $p$-group.
In particular this hypothesis holds if $p=2$ and for each $i$ the summands of $P_i$
have distinct exponents.
\end{Cor}
{\it Proof:} The first statement is clear from Corollary \ref{CorBAA}. Similarly the second follows 
from Proposition 4.5 of \cite{MPD} which implies each $Aut(P_i)$ is a $2$-group. \qed

\medskip

Another particularly tractable case is that of 
$E: N \to  P \to Q$,
with $N$ an Abelian, characteristic subgroup.
Since $Hom(Q,N)$ tends to be relatively large in this case, the quotient $Im(\rho)$ 
can be significantly 
smaller than $Aut(G)$. We examine this phenomenon in more detail.

Let $\{\Gamma^n(P)\}$ be the
{\it mod-p lower central series}, that is $\Gamma^0(P) = P$ and
$$
  \Gamma^n(P) = \langle (g_1, \dots, g_s)^{p^k} \mid sp^{k} > n \ \rangle, \ \ \ n\geq 1
$$
where $(g_1, \dots, g_s) = (g_1, (g_2, (\dots(g_{s-1}, g_s)\dots ))$ is the $s$-fold iterated
commutator. Then $\Gamma^1(P) = \Phi(P)$ is the Frattini subgroup and $V = 
P/\Gamma^1(P)$
is the largest elementary Abelian quotient of $P$. For some $n$, $\Gamma^n(P) = 1$.

\medskip

\begin{Ex}
\label{ExA}
\rm For $p$ an odd prime we consider the group 
\medskip

$P=\langle a,b,c \ | \ a^{p^2}= b^{p^2}= c^{p}= 1, c= (b,a)$, trivial higher commutators $ \rangle$ 
\medskip

We shall determine $Aut(P/\Gamma^p(P))$ using the extension
$$
  E: N = \Gamma^1(P)/ \Gamma^p(P) \to P/\Gamma^p(P) \to Q = P/\Gamma^1(P)
$$
where $N = \mathbf{Z}/p \langle a^p, b^p, c \rangle$ and 
$Q = \mathbf{Z}/p \langle \overline{a}, \overline{b} 
\rangle$ are elementary Abelian. This extension has trivial twisting and its 
extension cocycle is easily seen to be 
$[E] = (x\wedge y) \otimes c \in H^2(V,Q) = H^2(V) \otimes \mathbf{Z}/p \langle z \rangle$
where $x,y$ are dual to $a,b$ respectively. By inspection 
$$
Im(\rho) = \langle (A,B) \in GL_2(\mathbf{F}_p)\times GL_1(\mathbf{F}_p) \ | \ B= det(A) \rangle.
$$
Thus $Im(\rho) \cong GL_2(\mathbf{F}_p)$. We conclude 
$\mathbf{F}_pAut(P) \to \mathbf{F}_pGL_2(\mathbf{F}_p)$ is surjective with nilpotent kernel.
Finally we note $Hom(Q, N) = Mat_{3,2}(\mathbf{F}_p)$, thus $|Hom(Q,N)| = p^6$ and
$|Aut(P)| = p^6 |GL_2(\mathbf{F}_p)| = p^7(p^2 -1)(p-1)$.
\end{Ex}
\medskip

\begin{Ex}
\label{ExB}
\rm Let $U_5(\mathbf{F}_2) \leq GL_5(\mathbf{F}_2)$
be the unipotent subgroup of upper triangular matrices over $\mathbf{F}_2$.
These matrices necessarily have ones on the diagonal. Thus $U_5(\mathbf{F}_2)$ 
is generated by $x_{ij} = I_4 + e_{ij}$ where $1\leq i < j\leq 5$ and $e_{ij}$ 
is the standard elementary matrix with $1$ in the $ij$ position and zeros elsewhere.
Let $P = U_5(\mathbf{F}_2)/\Gamma^2U_5(\mathbf{F}_2)$. Then there is a central extension
$$
E: N = \mathbf{Z}/2\langle \overline{x}_{13}, \overline{x}_{24}, \overline{x}_{35} \rangle \to P \to  
    Q = \mathbf{Z}/2 \langle \overline{x}_{12}, \overline{x}_{23}, \overline{x}_{34},
    \overline{x}_{45} \rangle
$$
Since $N$ is the commutator subgroup of $P$, it is characteristic. The extension cocycle is
$$
[E] = y_{12}y_{23}\otimes \overline{x}_{13} + y_{23}y_{34} \otimes \overline{x}_{24} + y_{34}y_{45} 
\otimes \overline{x}_{35}
$$
where 
$$
    H^*(Q, N) = H^*(Q) \otimes N = 
\mathbf{Z}/2\langle y_{12}, y_{23}, y_{34}, y_{45} \rangle \otimes \langle \overline{x}_{13}, \overline{x}_{24}, \overline{x}_{35} \rangle
$$
and $y_{ij}$ is dual to $\overline{x}_{ij}$. Since the twisting is trivial, $C_{\chi} = Aut(Q)\times
Aut(N) = GL_4(\mathbf{F}_2)\times GL_3(\mathbf{F}_2)$.
Direct calculation shows 
$$
(Aut(Q) \times Aut(N))_{[E]} =  \langle A,B \ \ | \ \ A^4 = B^2 = 1, BAB = A^{-1} \rangle 
$$
the dihedral group of order $8$ generated by
\smallskip
$$
A = \left( \begin{matrix}    
0& 0& 0& 1\\
0& 0& 1& 0\\
0& 1& 0& 0\\
1& 0& 1& 0\\
\end{matrix} \right)
\times
\left ( \begin{matrix}
0& 0& 1\\
1& 1& 0\\
1& 0& 0\\
\end{matrix} \right)
$$

$$
B= \left( \begin{matrix}
1& 0& 0& 0\\
0& 1& 0& 0\\
0& 0& 1& 0\\
0& 1& 0& 1\\
\end{matrix} \right)
\times
\left ( \begin{matrix}
1& 0& 0\\
0& 1& 1\\
0& 0& 1\\
\end{matrix} \right)
$$
\medskip
\noindent Moreover $|Hom(Q,N)| = 2^{12}$, $Aut(P)$ is a $2$ group of
order $2^{15}$.
\end{Ex}
\bigskip

\section{An inductive procedure for determining $Aut(P)$}

Let $G=P$ be a $p$-group.
Since $\Gamma^i(P)$ is characteristic there is an induced homomorphism
$\rho_V: Aut(P) \to Aut(V)$ which factors as
$$
\rho_V:  Aut(P) \to \cdots \to Aut(P/\Gamma^{i+1}(P))\overset{\rho_{i}}\rightarrow Aut(P/\Gamma^i(P))
$$
$$
\to \cdots \overset{\rho_2}\rightarrow Aut(P/\Gamma^2(P))\overset{\rho_1}\rightarrow Aut(V)
$$
We shall describe an inductive procedure for lifting elements in the image of this map.

Consider the extensions
$$
    E_i: \Gamma^1(P)/\Gamma^{i}(P) \to P/\Gamma^{i}(P) \to V,  \ \ \ i\geq 2
$$
$$
{\widetilde E}_i: \Gamma^i(P)/\Gamma^{i+1}(P) \to \Gamma^1(P)/\Gamma^{i+1}(P) \to \Gamma^1(P)/\Gamma^i(P), \ \ \
i\geq 2.
$$
where $E_2$ and ${\widetilde E}_i$ have trivial twisting. In each case the kernel is a characteristic
subgroup.

Let $\sigma_1 \in Aut(V)$. By Theorem \ref{ThmWB},
$$
  \sigma_1 \in Im\{Aut_{\Gamma^1(P)/\Gamma^2(P)}(P/\Gamma^2(P)) \to Aut(V)\}
$$
if and only if there
exist $\tau_1 \in Aut(\Gamma^1(P)/\Gamma^2(P))$ such that
$(\sigma_1, \tau_1)$ fixes
$$
[E_2] \in H^2(V, \Gamma^1(P)/\Gamma^2(P)).
$$
Then there exists $\sigma_2  \in
Aut_{\Gamma^1(P)/\Gamma^2(P)}(P/\Gamma^2(P))$ lifting $\sigma_1$. Since
$$
Aut_{\Gamma^1(P)/\Gamma^2(P)}(P/\Gamma^2(P))= Aut(P/\Gamma^2(P))
$$
this completes the initial step.
Now suppose inductively that we have found elements $\sigma_i \in Aut(P/\Gamma^i(P))$,
$\tau_{i-1} \in Aut(\Gamma^1(P)/\Gamma^{i}(P))$ such
that $\sigma_i, \tau_{i-1}$ are  lifts of $\sigma_{1}, \tau_{i-2}$, respectively.
We need to find
$\tau_{i} \in Aut(\Gamma^1(P)/\Gamma^{i+1}(P))$ such that
$(\sigma_1, \tau_i)\in C_{\chi}$ fixes
$$
   [E_{i+1}] \in H^2_{\overline{\chi}}(V, Z( \Gamma^1(P)/\Gamma^{i+1}(P))).
$$
Then by Thm \ref{ThmWB}
there exists $\sigma_{i+1}\in Aut(P/\Gamma^i(P))$ lifting $\sigma_1$.

To find $\tau_{i}$ we apply the same technique to the extension $\widetilde{E}_{i}$
noting that the twisting is trivial so condition (1) is trivially satisfied.
Thus we must find ${\tau}'_i \in Aut(\Gamma^i(P)/\Gamma^{i+1}(P))$ such that
$(\tau_{i-1}, {\tau}'_{i})$ fixes
$$
[\widetilde{E}_{i}] \in H^2( \Gamma^1(P)/\Gamma^i(P), \Gamma^i(P)/\Gamma^{i+1}(P))
$$
since $ \Gamma^i(P)/\Gamma^{i+1}(P)$ is its own center.
Then applying Theorem \ref{ThmWB} again there exists $\tau_i \in Aut(\Gamma^1(P)/\Gamma^{i+1}(P))$
lifting $\tau_{i-1}$ as desired.

The induction terminates when $\Gamma^{i+1}(P) = 1$.

The procedure described in this section is demonstrated in the examples in Section 7.

\bigskip

\section{Extensions of elementary Abelian groups}

We consider extensions  $ E: N \to G \to V$ where $N,V$ are elementary
Abelian $p$-groups with $N$ central.
In this case the twisting is trivial, $\chi = id$ and $C_{\chi} = Aut(N) \times Aut(V)$.
Our aim is to use Corollary \ref{CorABB} to compute $|Im(\rho)|$. 
Let $n= dim_{\mathbf{F}_p}(N)$. Then the extension cocycle 
$[E] \in H^2(V;N)$ has the form
$[E] = (X_1, X_2, \dots ,X_n)$ where  $X_i \in H^2(V; \mathbf{F}_p)$.
We recall that $Aut(V)$ acts diagonally on $[E]$, $\sigma [E] = (\sigma X_1, \sigma X_2,
\dots, \sigma X_n)$ for $\sigma \in Aut(V)$. Thus the isotropy subgroup 
$$
Aut(V)_{[E]} = Aut(V)_{X_1} \cap \cdot \cdot \cdot \cap Aut(V)_{X_n}\ \ \ \ \ \ \ \ \ \ \ \ \ \ \ \ (3)
$$
The action of $Aut(N)$ on $[E]$ is induced from that on $N$.
\medskip

\subsection{ Quadratic Forms}

\smallskip

At this point we restrict our attention to the case $p=2$. Let $m = dim(V)$ then each $X_i$
is a quadratic form in $x_1, x_2, \dots , x_m$ the generators of 
$H^*(V; \mathbf{F}_2) = \mathbf{F}_2[x_1, x_2, \dots , x_m]$, $|x_i| = 1$.

We recall some classical facts about quadratic forms $Q: V \to \mathbf{F}_2$, \cite{Br},
\cite{Di}, \cite{Dd}.
The defining property is that
the associated form $B(x,y) = Q(x+y) + Q(x) +Q(y)$ is alternate bilinear.

The {\it bilinear radical} of $B$,           
$$
bilrad(V,B) := \{ x \in V | B(x,y) = 0,  \forall y \in V \}
$$
As usual, $B$ is called {\it non-degenerate} if $bilrad(V,B) = 0$, i.e. its matrix is non-singular.
The {\it radical} of $Q$, 
$$
Rad(V,Q) := \{x \in bilrad(V,B) | Q(x) =0 \}
$$
$Q$ is said to be {\it non-degenerate} if $Rad(V,Q) = 0$. 

By a theorem of Dickson  \cite{Di} (Section 199) a (non-zero) quadratic form over $\mathbf{F}_2$
in $m$ variables which is not equivalent (by a change of basis) to one in fewer variables
must be equivalent to one of the following standard non-degenerate quadratic forms
$$
\Phi_{m}^+ =   x_1 x_2 + \cdot  \cdot  \cdot + x_{m-1} x_{m}, \ \ \ \ \ \ \ \ \ \ \ \ \ \ \ \ \ \ \ \ \ \ \ \ \ \ \ \ \ \ \ \ \ \ \ \ \ \ m  \ \ even
$$
$$
\Phi_{m}^- =  x_1 x_2 + \cdot  \cdot  \cdot + x_{m-3} x_{m-2} + {x_{m-1}}^2 + x_{m-1} x_{m} + {x_{m}}^2, \ \ m \ \ even
$$
$$
\Phi_m =  {x_1}^2 + x_2 x_3 + \cdot  \cdot  \cdot + x_{m-1}x_{m},  \ \ \ \ \ \ \ \ \ \ \ \ \ \ \ \ \ \ \ \ \ \ \ \ \ \ \ \ \ \ \ \ m \ \ odd
$$

If $bilrad(V,B) = 0$,
then $m=2r$ is even and one can define the {\it Arf invariant} of $Q$ with respect
to a symplectic basis $\{u_1, v_1,\dots, u_{k}, v_{k}\}$ by
$$
     Arf(Q) = \sum_{i=1}^k Q(u_i)Q(v_i) \in \mathbf{Z}/2
$$
This is invariant of the choice of symplectic basis and determines $Q$ up to equivalence.
It is convenient to write ${\mathbf Z}/2 = \langle \pm 1 \rangle$ multiplicatively. Then with this notation,
W. Browder has shown that $Arf(Q) = 1, -1$ if and only if $Q$ sends the majority of elements
of $V$ to $1, -1$ respectively \cite{Br}. 
For $m$ even one finds $Arf(\Phi_{m}^+) = 1$, $Arf(\Phi_{m}^-) = -1$.

The Arf invariant can be extended to the 
$m$ odd case (where $B$ is degenerate) as follows. It is easy to see that $\Phi_m$ 
sends the same number of elements to $1$ and $-1$ thus one can define $Arf(\Phi_m) = 0$.
It is clear that Browder's definition (also known as the ``democratic invariant") 
is invariant under any basis change. This leads to the following classification theorem.

\begin{Thm}
\label{ThmWY}\cite{Br},Theorem III.1.14

\noindent A quadratic form $Q:V \to \mathbf{Z}/2$ is determined up to equivalence by the triple 
$(dim(V), dim(bilrad(V,B)), Arf(Q))$

\end{Thm}

\medskip

The action of $Aut(N) = GL(n,\mathbf{F}_2)$ on an $n$-tuple $(X_1, X_2,\dots ,X_n)$ of quadratic forms 
is linear and thus involves the sum of forms. Unfortunately it is impossible, in general, to determine 
the sum from the Arf invariant. For example if $X_1 = x^2$ and $X_2 = xy +y^2$ then $Arf(X_1 + X_2)
\neq Arf(X_1) + Arf(X_2)$. However, on a direct sum of vector spaces then it follows easily from the 
definition that 
$$
    Arf(X_1 \oplus X_2) = Arf (X_1) + Arf(X_2).
$$  

\bigskip

For the rest of this subsection we restrict attention to $m$ even and 
consider the special case $[E] =  (X_1, X_2, \dots ,X_n)$ where the $X_i$ are in standard form.
We write $X_i = X$ or $Y$ depending on whether the
Arf invariant is $1$ or $-1$. In the following theorem we use $\sim$ to 
denote conjugacy.
\medskip

\begin{Thm}
\label{ThmXY} Suppose $m=2r$ and $[E] = (X_1, X_2, \dots ,X_n)$ with $X_1 =\cdots = X_k = X$,
$X_{k+1} = \cdots = X_n = Y$. Then $\Omega([E]) = [E]$ and 
$|Im(\rho)| = |Aut(V)_{[E]}|\cdot |Aut(N)_{[E]}|$.
Furthermore

1) If $k=n$ then
$$
    Aut(V)_{[E]} = O_{m}^+(\mathbf{F}_2)
$$
the orthogonal group of order $2(2^r - 1)\prod_{i=1}^{r -1}(2^{2i} -1)2^{2i}$ of matrices 
preserving the form $X$.
$$
      Aut(N)_{[E]}\sim
\left ( \begin{matrix}
1 & 0\\
*& GL_{n-1}(\mathbf{F}_2)\\
\end{matrix} \right)
$$
of order $2^{n-1}\prod_{i=1}^{n-1}(2^i - 1)2^{i-1}$.
$Im(\rho)$ is a $2$-group if and only if $m,n < 3$.

2) If $k = 0$ then
$$
    Aut(V)_{[E]} = O_{m}^-(\mathbf{F}_2)
$$
the orthogonal group of order $2(2^n + 1)\prod_{i=1}^{n -1}(2^{2i} -1)2^{2i}$ of matrices 
preserving the form $Y$ and
$$
     Aut(N)_{[E]} \sim
\left ( \begin{matrix}
1 & 0\\
*& GL_{n-1}(\mathbf{F}_2)\\
\end{matrix} \right).
$$
$Im(\rho)$ is a $2$-group if and only if $m<2$ and $n < 3$.

3) If $1 < k <n$ then
$$
     Aut(V)_{[E]} =  O_{m}^+(\mathbf{F}_2) \cap  O_{m}^-(\mathbf{F}_2)
$$ 
and
$$
    Aut(N)_{[E]} \sim
\left ( \begin{matrix}
I_2& 0\\
* & GL_{n-2}(\mathbf{F}_2)\\
\end{matrix} \right)
$$
of order $2^{2(n-2)}\prod_{i=1}^{n-2}(2^i - 1)2^{i-1}$.
$Im(\rho)$ is a $2$-group if and only if $m,n <4$.

\end{Thm}

{\it{Proof:}}
The calculation of $Aut(V)_{[E]}$ follows from (3). The intersection orbit group
$\Omega([E])=
(Aut(V)\cdot [E]) \cap (Aut(N)\cdot [E]) = \{ [E] \}$. To see this suppose
$\sigma [E] = [E]{\tau}^{-1}$ differs from $[E] = (X_1,X_2,\dots,X_n)$ in the $i$-th coordinate $X_i = X$, say.
Now $\sigma X = ([E]{\tau}^-1)_i = aX + bY$, $a,b \in \mathbf{F}_2$.
Thus $\sigma X = X + Y$ since $\sigma X \neq X$ and $Y$ has Arf invariant $-1$.
However $\sigma X = X + Y = x_{m-1}^2 + x_{m}^2 = (x_{m-1} + x_{m})^2$
contradicting the fact that $X$ is not equivalent to a quadratic form in fewer than
$m$ variables. Similarly if $X_i = Y$.

Since $X$ and $Y$ are linearly independent polynomials, $(X,\dots,X,Y,\dots,Y)$
is equivalent to $(X,Y,0,\dots,0)$ by a change of basis. Thus the descriptions of  $Aut(N)_{[E]}$
follow immediately. \qed

\bigskip

\subsection{ Non-standard Forms}

\smallskip

We now turn to the more general case where the forms $X_i$ of $[E] = (X_1, X_2, \dots ,X_n)$ are not 
in standard form.
One can still determine  $|Im(\rho)|$; 
we illustrate this by considering the case $n =2$, $m=3$.
Thus we consider pairs of quadratic forms $(X,Y)$ in variables $x,y,z$ each equivalent to
(but not necessarily equal to) one of the standard four forms:
$$
 x^2 + yz, \ \ \  xy, \ \ \  x^2 + xy + y^2, \ \ \ x^2.
$$
The respective isotropy groups, as subgroups of $GL_3(\mathbf{F}_2)$, are
$$
Aut(V)_{x^2 + yz} = O_3(\mathbf{F}_2) \cong GL_2(\mathbf{F}_2)
$$ 
non-Abelian of order $6$;
$$
Aut(V)_{xy} = 
\left( \begin{matrix}
\Sigma_2& 0\\
*& 1\\
\end{matrix} \right)
$$
elementary Abelian of rank $3$;
$$
Aut(V)_{x^2 +xy + y^2} = 
\left( \begin{matrix}
GL_2 & 0\\
*& 1\\
\end{matrix} \right)
$$
of order $24$;
$$
Aut(V)_{x^2} = 
\left( \begin{matrix}
1 & 0\\
* & GL_2\\
\end{matrix} \right)
$$
of order $24$.
\bigskip

\subsubsection{ \it Simultaneous Equivalence}

First we consider the case where $X$ and $Y$ are simultaneously equivalent to a standard
form, i.e., there is an invertible linear transformation $A$ of $V$
such that $A^{-1}XA$ and $A^{-1}YA$ are each in standard form.
\medskip

\noindent Case 1, $X = Y$: 

Then $Aut(V)_{(X,X)} = Aut(V)_{X}$. It is also easy to see
$Aut(N)_{(X,X)} = \Sigma_2$ and the intersection of orbits $\Omega([E]) = \{(X,X)\}$ in 
this case.
\medskip

a) $X$ is equivalent to $xy:$ As above the isotropy subgroup $Aut(V)_{(X, X)}$ is
elementary Abelian of rank $3$. The intersection of the orbits
is just $(X,X)$ so the order of $Im(\rho) = 8\cdot 2 \cdot 1 = 16$.
\medskip

b) $X$ is equivalent to $x^2 + xy + y^2:$ $Aut(V)_{(X,X)}$ is of order $24$.
Thus $|Im(\rho)| = 24\cdot 2 \cdot1 = 48$.
\medskip

c) $X$ is equivalent to $x^2 + yz:$ $Aut(V)_{(X,X)} = O_{3}(F_2)$. 
Thus $|Im(\rho)| = 6 \cdot 2 \cdot 1 = 12$.
\medskip

d) $X$ is equivalent to $x^2$: As above $Aut(V)_{(X,X)}$ is of order $24$. Thus  
$|Im(\rho)| = 24 \cdot 2 \cdot 1 = 48$.
\bigskip

\noindent Case 2,  $X \neq Y$: 

In this case $Aut(N)_{(X, Y)} = 1$. 
Further $Aut(V)_{(X,Y)} =
Aut(V)_{X} \cap Aut(V)_{Y} = Aut(V)_{(Y,X)}$. There are several possibilities; we give only four
in detail since the rest follow the same general pattern. By considering the Arf invariant we see that
$\Omega([E]) = 1$ except in the third example.

\medskip

1) $(X,Y)$ equivalent to $(x^2 + yz, x^2 + xy + y^2)$:  Then $Aut(V)_{X} \cap Aut(V)_{Y}$ is
$$
O_3(\mathbf{F}_2) \cap
\left( \begin{matrix}
GL_2& 0\\
*& 1\\
\end{matrix} \right)
= \mathbf{Z}/2 \Bigg\langle
\left( \begin{matrix}
1& 1 &0\\
0& 1 &0\\
0 & 1 & 1\\
\end{matrix} \right)
\Bigg\rangle.
$$
The intersection of orbits is $\{(X,Y)\}$ hence  $|Im(\rho)| = 2 \cdot 1 \cdot 1 = 2$.
\medskip

2) $(X,Y)$ equivalent to $(xy, x^2 +xy +y^2)$: Then $Aut(V)_{X}\cap Aut(V)_{Y}$ is
$$
\left( \begin{matrix}
\Sigma_2 & 0\\
* & 1\\
\end{matrix} \right)
\cap
\left( \begin{matrix}
GL_2& 0\\
*& 1\\
\end{matrix} \right)
=
\left( \begin{matrix}
\Sigma_2 & 0\\
* & 1\\
\end{matrix} \right)
$$
i.e., dihedral of order $8$. The intersection of orbits is $\{(X,Y)\}$ thus  $|Im(\rho)| = 8 \cdot 1 \cdot 1 = 8$.

Similarly for $(xy, y^2 + yz + z^2)$ and $(xz, x^2 + xy + y^2)$.
\medskip

3) $(X,Y)$ equivalent to $(xy, x^2 + yz)$: Then $Aut(V)_X \cap Aut(V)_Y$ is trivial
and $\Omega([E]) = \{(X,Y), (X+Y,Y) \} = \mathbf{Z}/2$. Thus  $|Im(\rho)| = 2$.
\medskip

4) $(X,Y)$ equivalent to $(xy, yz)$: Then $Aut(V)_{X} \cap Aut(V)_{Y}$ is trivial.
Direct calculation shows the intersection of the orbits
has order $6$. Thus  $|Im(\rho)| = 1 \cdot 1 \cdot 6 = 6$.

Similarly for $(xz, yz)$, $(xz, xy)$, and $(yz, xz)$.
\medskip

\medskip

\subsubsection{\it Non-simultaneous Equivalence}

By this we mean $X$ and $Y$ are not simultaneously equivalent to a pair of standard forms.
Since $X \neq Y$, $Aut(N)_{(X,Y)} = 1$ in all cases.

To illustrate this phenomena the following table gives a complete computation of $|Im(\rho)|$ 
in case $X$ and $Y$ are equivalent (non-simultaneously) to $x^2 + yz$. Then $X$ may be assumed
to be $x^2 + yz$ and we list only the relevant $Y$'s and the corresponding values of $|Im(\rho)|
= 2, 3, 4$.
\bigskip

\centerline{\bf{$|Im(\rho)|$}}
\medskip
\begin{tabular}{|c|c|c|} \hline
2 & 3 & 4 \\ \hline
$x^2 +xz + yz$ & $xy + xz + y^2$ & $x^2 + yz + z^2$ \\ 
$x^2 + xy + xz + z^2$ & $x^2 + xy + xz + y^2 + yz$ & $xy + xz + y^2 + yz + z^2$ \\ 
$x^2 + xz + yz + z^2$ & $xz + y^2 + yz$ & $ xy + xz + yz$ \\ 
$x^2 + xy + y^2 + yz$ & $xy + xz + z^2$ & $xz + y^2 + z^2$  \\ 
$x^2 + xy + yz$ & $ x^2 + xy + y^2 + z^2$ & $x^2 + y^2 + yz$ \\ 
$x^2 + xy + xz + y^2$ & $x^2 + xz + y^2 + z^2$ & $xy + z^2$ \\ 
     & $x^2 + xy + z^2$ & $xy + y^2 + z^2$ \\ 
  &    $xz + y^2 + yz + z^2$ & $ x^2 + y^2 + yz + z^2$ \\ 
  &    $x^2 + xz + y^2$ & $xz + y^2 $ \\
  &    $x^2 + xy + xz + yz + z^2$ & \\ 
  &    $xy + y^2 + yz + z^2$ & \\ 
  &    $   xy+ yz + z^2$   & \\ \hline
\end{tabular}

\bigskip

In more detail,
if $|Im(\rho)| = 2$, then  $Aut(V)_{(X,Y)} = 1$ and 
$\Omega([E]) = \{(X,Y), (Y,X)\} = \mathbf{Z}/2$
If $|Im(\rho)| = 3$, then $ |Aut(V)_{(X,Y)}| = 1$ and $\Omega([E]) = 
\{(X, Y), (X+Y, X), (Y, X+Y) \} = \mathbf{Z}/3$. If $|Im(\rho)| = 4$ then $|Aut(V)_{(X,Y)}| = 2$
and $\Omega([E]) = \{(X,Y), (Y,X) \} = \mathbf{Z}/2$.

\bigskip

As a final example we consider
\medskip

\noindent $X = xy + y^2$, $Y = x^2 + y^2 + yz + z^2$.
To analyze $Aut(V)_{(X,Y)}$ we separately reduce $X$, $Y$ to standard forms $f_1, f_2$
respectively
$$
\sigma_1 X = f_1, \ \ \ \sigma_2 Y = f_2
$$
with $\sigma_i \in Aut(V)$. Let $\sigma = \sigma_2{\sigma_1}^{-1}$, then ${\sigma}^{-1}f_2
= \sigma_1 Y$ hence
$$
      Aut(V)_{(f_1, \sigma_1 Y)} = Aut(V)_{f_1} \cap {\sigma}^{-1}[Aut(V)_{f_2}] \sigma
$$
This determines $Aut(V)_{(X,Y)}$ up to conjugacy since
$Aut(V)_{(X,Y)} = \\{\sigma_1}^{-1}[ Aut(V)_{(f_1, \sigma_1 Y)} ]\sigma_1$. In this example
$f_1 = xy$, $f_2 = x^2 + yz$ with
$\sigma_1: y \mapsto x+ y$, $\sigma_2 = x \mapsto x + y +z $. Then
$$
Aut(V)_{(f_1, \sigma_1 Y)}=Aut_{xy}\cap {\sigma}^{-1}[Aut(V)_{x^2 + yz}] \sigma = 
\mathbf{Z}/2
\Bigg\langle \left( \begin{matrix}
1& 0 &0\\
0& 1 &0\\
0 & 1 & 1\\
\end{matrix} \right)\Bigg\rangle
$$
Further $\Omega([E]) = \{(X,Y),(X,X+Y)\} = \mathbf{Z}/2$ hence $|Im(\rho)| = 2 \cdot 1 \cdot 2 = 4$.
\medskip


\bigskip

\section{Applications}

\noindent{\bf{1.}}
Let $P$ denote the extraspecial group of order $|P| = p^{2n+1}$ and exponent $p>2$
defined by the central extension
$$
   E: Z \to P \to V
$$
where $Z = \Phi(P) = \mathbf{Z}/p$, $V = P/\Phi(P) = (\mathbf{Z}/p)^{2n}$.
The twisting $\chi$ is trivial thus $C_{\chi} =  Aut(V) \times Aut(Z)$.
Complete information about the automorphism group of $P$ as well as all other extensions
of elementary Abelian $p$-groups by $\mathbf{Z}/p$ is known \cite{Win},\cite {D}.
Here we apply our results to obtain a quick derivation of 
$Im\{Aut(P) \overset{\rho}\rightarrow Aut(V) \times Aut(Z)\}$.

$P$ is generated by elements $x_1, x_2, \dots , x_{2n}, \zeta$ of order $p$ satisfying
$$
  [\zeta, x_i] = 1
$$
$$
   [x_{2i-1}, x_{2i}] = \zeta
$$
$$
  [x_{2i-1}, x_j ] = 1, \ \ j \neq 2i
$$
$$
  [x_{2i}, x_j ] = 1, \ \ j \neq 2i-1
$$
Thus $\Phi(P) = Z(P) = \langle \zeta \rangle = \mathbf{Z}/p$
and $V = \langle \overline{x}_1, \overline{x}_2, \dots , \overline{x}_{2n} \rangle$.
It is immediate from these relations that the extension cocycle is
$$
[E] = B \otimes \zeta
$$
where $B =  y_1y_2 + \cdots + y_{2n-1}y_{2n}$,
$y_i \in H^1(V)$ is dual to $\overline{x}_i$.
Since $p$ is an odd prime, ${y_i}y_j = -y_j{y_i}$. Thus
$B$ is exactly the skew-symmetric form for the symplectic group
$Sp(n, \mathbf{F}_p)$. We shall need a slightly more general version $GSp(2n, \mathbf{F}_p)$,
\cite{Dd}, the transformations which fix $B$ up to a scalar. It is easy to see that
$GSp(2n, \mathbf{F}_p) = \langle \gamma \rangle \rtimes  Sp(2n, \mathbf{F}_p)$ where
$\gamma$ is the linear transformation
$$
    \overline{x}_{2i-1} \mapsto k \overline{x}_{2i-1},\ \ \ \ \  \overline{x}_{2i} \mapsto  \overline{x}_{2i}
$$
and $k$ is a generator of ${\mathbf{F}_p}^*$.

Now considering $(\sigma, \tau) \in Aut(Z) \times Aut(V)$ acting on $E$
we find
$$
(\sigma, \tau)(E) = \lbrack \sum_{i=1}^{n}(\sigma)(y_{2i-1})\sigma(y_{2i})\rbrack\otimes \tau(\zeta)
$$
Thus $ (\sigma, \tau)(E) = E$ if and only if $\sigma \in GSp(2n, \mathbf{F}_p)$ and $\tau$ is multiplication
by $det(\sigma)^{-1}$. We conclude $Im\{Aut(P) \to Aut(V) \times Aut(Z)\} \cong
GSp(2n, \mathbf{F}_p)$.
\bigskip

\noindent{\bf{2.}}
Let $W(n)$ be the {\it{universal W-group}} \cite{AKM} on $n$ generators defined as the 
central extension
$$
   1 \to N = \Phi(W(n)) \to W(n) \to Q = (\mathbf{Z}/2)^n \to 1
$$
where $N = (\mathbf{Z}/2)^{n + \binom n2}$. These groups arise, for instance, in the study of the 
cohomology of Galois groups.
The extension class $[E] = (X_1, X_2, \dots, X_{n + \binom n2})$
where the $\{X_i \}$ form an ordered basis for $H^2(Q)= S^2(Q^*)$, the second symmetric power of
the dual of $Q$. Any order will do. The twisting is trivial, hence $C_{\chi} = GL(Q) \times GL(N)$.
Now $GL(Q)$ induces linear isomorphisms on $H^*(Q)$ which are determined by their values on squares
in $H^2(Q)$.
Thus given $\sigma \in GL(Q)$ we can find a
$\tau \in GL(N)$ such that $(\sigma, \tau)([E]) = (\sigma [E]){\tau}^{-1} = [E]$. 
Thus we see that $Im(\rho) \leq  GL(Q) \times GL(N)$ by an injection
which projects to an isomorphism on the first factor, i.e. $Im(\rho)\cong GL(Q)$.

\bigskip

\noindent{\bf{3.}}
We consider $Aut(P)$ for the unipotent group 
$P = U_4(\mathbf{F}_2) \subset GL_4(\mathbf{F}_2)$
of upper triangular matrices over the finite field $\mathbf{F}_2$. These 
matrices necessarily have ones on the diagonal. Thus $P$ is generated by 
$\{x_{12}, x_{23}, x_{34}\}$. Then $\Gamma^3(P) = 1$, $\Gamma^2(P) = 
\mathbf{Z}/2 \langle x_{14}\rangle = ZP$, and $\Gamma^1(P) = 
\mathbf{Z}/2 \langle x_{13}, x_{24}, x_{14} \rangle$. The corresponding extensions 
have
$$
E_2:  \Gamma^1(P)/\Gamma^2(P) \overset{i_2}\rightarrow  P/\Gamma^2(P) \to V
$$
\smallskip
\noindent where $V = P/\Gamma^1(P) = 
\mathbf{Z}/2 \langle \overline{x}_{12}, \overline{x}_{23}, \overline{x}_{34} \rangle$, 
$ \Gamma^1 (P) / \Gamma^2 (P) 
= \mathbf{Z}/2 \langle \overline{x}_{13}, \overline{x}_{24} \rangle $
and
$$
E_3:  \Gamma^1(P) \overset{i_3}\rightarrow P \to V
$$
The twisting $\chi$ is trivial for $E_2$ but not for $E_3$.
In this application we shall not need the auxiliary extensions $[\widetilde{E}_i]$.

\medskip
\noindent{\bf Claim:} $Im \{ Aut(P) \to Aut(V) \} \cong \Sigma_3$
\medskip

\noindent{\it Proof:} (Sketch) From the commutator relations 
$[x_{ij}, x_{jk}] = x_{ik}$ for $i< j <k$,
one sees that the extension class $[E_2] \in H^2(V; \Gamma^1(P)/\Gamma^2(P))$
is $xy\otimes \overline{x}_{13} + yz\otimes \overline{x}_{24}$ where
$$
   H^*(V,  \Gamma^1(P)/\Gamma^2(P)) = H^*(V) \otimes \mathbf{Z}/2 \langle {x}_{13}, {x}_{24} \rangle 
$$
$$
= \mathbf{Z}/2[x,y,z] \otimes 
\mathbf{Z}/2 \langle {x}_{13}, {x}_{24} \rangle
$$   
Here $x,y,z\in H^1(V)$ are classes dual to 
$ \overline{x}_{12}, \overline{x}_{23}, \overline{x}_{34}$ respectively.
Direct calculation shows the $[E_2]$ is fixed the subgroup
generated by the involution
defined by
$$ 
\sigma:\   \overline{x}_{12} \mapsto \overline{x}_{34}, 
\ \  \overline{x}_{23} \mapsto \overline{x}_{23}, 
\ \ \tau: \   x_{13} \mapsto x_{24}
$$ 
and the map of order three 
$$
\sigma': \ \overline{x}_{12} \mapsto \overline{x}_{12} + \overline{x}_{34}, 
\ \ \overline{x}_{23} \mapsto \overline{x}_{23}, 
\ \  \overline{x}_{34} \mapsto \overline{x}_{12},
$$
$$
\tau':\  x_{13} \mapsto x_{24}, \ \ x_{24} \mapsto x_{13} + x_{24}
$$
One finds $\langle (\sigma, \tau), (\sigma', \tau') \rangle \cong \Sigma_3$.

Turning to the extension $E_3$, we can extend $\tau$, $\tau'$ to $\Gamma^1(P)$
by letting them act identically on $x_{14}$. Then a simple calculation shows
$\langle (\sigma, \tau), (\sigma', \tau') \rangle \leq C_{\chi}$. 
We note that $[E_3] = [E_2]$ considered as an
element of $H^2_{\overline{\chi}}(V;\Gamma^1(P))$. 
This follows from the definition of the cocycle 
$$
f: V \times V \to \Gamma^1(P)/\Gamma^2(P)
$$
for $[E_2]$. Recall that given a set theoretic section $s: V \to P/\Gamma^2(P)$ 
then
$$
s(a)s(b) = i_2(f(a,b))s(ab)
$$

If we use one of the sections not involving the center, 
then it lifts to  a section  $\tilde{s}:V \to P$
and the corresponding cocycle $\tilde{f}:V \times V \to \Gamma^1(P)$ is a lift of $f$.
(For example using the ordered basis $(\overline{x}_{12}, \overline{x}_{23}, \overline{x}_{34})$ 
for $V$ let $\tilde{s}(\overline{x}_{12}) = 
x_{12},\ \tilde{s}(\overline{x}_{23})=x_{23},\ \tilde{s}(\overline{x}_{34}) = x_{34}$.
For products $w \in V$ let
$\tilde{s}(w) = w'$
where $w'$ is ordered lexicographically. Thus if 
$w = \overline{x}_{34}\overline{x}_{23}\overline{x}_{12}$ then $w' = x_{12}x_{23}x_{34}$.
Hence 
$$
    \tilde{s}(\overline{x}_{34}\overline{x}_{23})\tilde{s}(\overline{x}_{12}) =
i_3(\tilde{f}(\overline{x}_{34}\overline{x}_{23}, \overline{x}_{12})) \tilde{s}(\overline{x}_{34}\overline{x}_{23}\overline{x}_{12})
$$
implies $\tilde{f}(\overline{x}_{34}\overline{x}_{23}, \overline{x}_{12})= x_{13}x_{24}$
In general
we observe that no two-fold products among the elements of 
$\tilde{s}(V)$ involve the center i.e., two-fold products among the elements 
$x_{12},\ x_{23},\ x_{34}$ 
are on the first and second diagonals, not in the center.)

This means $(\sigma, \tau), (\sigma', \tau') $ fix the extension class $[E_3]$ 
and hence define automorphisms of $P$ not just $P/\Gamma^1(P)$.

\qed

\medskip

\noindent{\bf 4.} Let $p$ be an odd prime and let
\medskip

$P=\langle a,b,c, d \ | \ a^{p}=b^{p}=c^{p}=d^p =1, c=(b,a), d = (c,a),$
 
$\ \ \ \ \ \ \ \ \ \ \ \ \ \ $ other commutators trivial$  \rangle$
\medskip

We shall study $Aut(P)$ for $p=5$, the  smallest prime for which $P$ is regular.
\noindent We consider the extensions  of $V = P/\Gamma^1(P) = \langle \overline{a}, \overline{b} \rangle$ 
$$
   [E_2]:  \ \ \  N_2 := \Gamma^1(P)/\Gamma^2(P)\to P/\Gamma^2(P) \to V
$$
$$
  [E_3]: \ \ \ \  \ \ \ \ N_3 := \Gamma^1(P) \to P \to V
$$
Since $\Gamma^{3}(P) = 1$, we also have
$$
[\widetilde{E}_2]:  \ \ \ \widetilde{N}_{3} := \Gamma^2(P) \to \Gamma^1(P) \to N_2 
$$
Each $N$ is Abelian:
$$
 N_2 = \langle \overline{c} \rangle, \ \ \ N_{3} = \langle c,d \rangle, \ \ \
\widetilde{N}_{3} =  \langle d \rangle
$$ 

Noting that $E_3$ has non-trivial twisting we apply the algorithm of Section 4 to study $Aut(P)$.

\medskip

\noindent{\bf Claim:} The semi-simple quotient of $\mathbf{F}_p(Aut(P))$
is $ \mathbf{F}_p(\mathbf{Z}/4)^2 $.
\medskip

\noindent{Proof:} (sketch)

Extension $E_2$ has trivial twisting and is quite similar to that of Example 1.
Using the same notation and arguing analogously we can determine the extension class 
$[E_2] = xy\otimes z \in H^2(V, N_2)$.
Hence 
$$
C_{[E_2]} = \langle (A,B)\in GL_2(\mathbf{F}_p)\times 
GL_1(\mathbf{F}_p) \ \ | \ \ B = det(A) \rangle 
$$
Extension $\widetilde{E}_2$ splits, $\Gamma^1(P) = \widetilde{N}_2 \times N_2 = \langle c, d \rangle$.
Thus every $\tau \in Aut(N_2)$ lifts to $Aut(\Gamma^1(P))$ and we proceed to study extension
$E_3$. 

First we determine 
$$
C_{\chi} = \langle (\sigma, \tau) \in Aut(V)\times Aut(N_3) \ \ 
| \ \ \chi \sigma = c_{\overline{\tau}} \chi \rangle
$$
The action of $V$ on $N_2$ is given by $\chi(\overline{x}): c \mapsto c + d, d \mapsto d$,
$\chi(\overline{y}) = id$. Hence $ker(\chi) = \langle \overline{y} \rangle$. Since 
$\sigma: ker(\chi) \to ker(\chi)$, it must have the form
$$ 
\sigma = \left( \begin{matrix} 
k& 0\\
m& n\\
\end{matrix} \right)
\in GL_2(\mathbf{F}_p)
$$
Now solving $\chi \sigma = c_{\overline{\tau}} \chi$ where
$$
\tau= \left( \begin{matrix} 
s& t\\
u& v\\
\end{matrix} \right)
\in GL_2(\mathbf{F}_p)
$$
we find
$$
\left( \begin{matrix} 
1& 0\\
k& 1\\
\end{matrix} \right)
=
\left( \begin{matrix} 
1+tv/\Delta& t^2/\Delta\\
v^2/\Delta& 1-tv/\Delta\\
\end{matrix} \right)
$$
where $\Delta = det(\tau)$.  Thus $t = 0$ and $k = v/s$.  Hence 
$$
C_{\chi} = \left\{\left( \begin{matrix}
k& 0\\
m & n\\
\end{matrix} \right)
\times
\left( \begin{matrix}
s& 0\\
u & v\\
\end{matrix} \right)  \ \bigg\vert\  \ k= v/s \right\}
$$
To determine $Im(\rho) = C_({\chi})_{[E_3]}$ it is easier 
in this case to observe that $|Aut(P)| = 4^2 5^5$ (using Magma) and argue directly instead of
using the extension class. First $ker(\rho)$
contains a subgroup of order $5^3$ generated by the inner automorphisms of order 
$5$, $\{ c_b, c_c \} \leq ker(\rho)$ together with the automorphism 
$\phi : \phi(b) = bd^{-1}, \phi = id$ on $a,c,d$. Next we show that $Im(\rho)$
contains a subgroup of order $4^2 5^2$ which must equal $Im(\rho)$.


To complete the proof of the claim we note that the subgroup of order $5^2$ is normal.\qed



\end{document}